\documentclass{article}
\usepackage[latin1]{inputenc}
\usepackage{amsmath, amssymb, amsthm}
\usepackage{amsfonts}

\def \({\left(}
\def \){\right)}

\def \C{\mathbb C}

\theoremstyle{definition}
\newtheorem{defi}{Definition}
\newtheorem{lemma}[defi]{Lemma}

\newtheorem{prop}[defi]{Proposition}

\title{A counterexample to the existence of a Poisson structure on a twisted group algebra}
\author{Eliana Zoque}
\date{}

\begin{document}
\maketitle
\begin{abstract}
Crawley-Boevey \cite{CB} introduced the definition of a noncommutative Poisson structure on an associative algebra $A$ that extends the notion of the usual Poisson bracket. Let $(V,\omega)$ be a symplectic manifold and $G$ be a finite group of symplectimorphisms of $V$. Consider the twisted group algebra $A=\C[V]\#G$.
We produce a counterexample to prove that it is not always possible to define a noncommutative poisson structure on $\C[V]\#G$ that extends
the Poisson bracket on $\C[V]^G.$
\end{abstract}

\section{Introduction}

Crawley-Boevey \cite{CB} defined a noncommutative Poisson structure
on an associative algebra $A$ over a ring $K$ as a Lie
bracket $\langle -,- \rangle$ on $A/[A,A]$ such that for each $a\in
A$ the map $\langle\overline a,-\rangle:a/[A, A]\to A/[A,A]$ is
induced by a derivation $d_a:A\to A$; i.e. $\langle \overline{a},,
\overline{b}\rangle=\overline{d_a(b)}$ where the map
$a\mapsto\overline{a}$ is the projection $A\to A/[A,A].$ When $A$ is
commutative a noncommutative Poisson structure is the same as a
Poisson bracket.

Let $(V,\omega)$ be a symplectic manifold, with the usual Poisson bracket
$\{-,-\}$ on $\C[V]$. Let $G$ be a finite group of
symplectimorphisms of $V$. Consider the twisted group algebra $A=\C[V]\#G$.
The algebra of $G$-invariant polymonials $\C[V]^G$ is contained in
$A/[A,A].$ We produce a counterexample to prove that it is not always possible to define a noncommutative poisson structure on $\C[V]\#G$ that extends
the Poisson bracket on $\C[V]^G.$

\section{Twisted group algebra and derivations}

From now on, let $A=\C[V]\#G$.($\C$ can be replaced by any field of characteristic 0.)
We use the symbol $^g\psi$ to denote the left action of $g\in G$ on $\psi\in \C[V]$. For every $g\in G$ we denote $(\cdot)_g$ the projection $A\to \C[V]$ into the $g$-part, i.e, $(\psi h)_g=\psi\delta_{g,h}$ if $\psi\in\C[V], h\in G.$ Let $G=C_0\cup C_1\cup\cdots$ be the conjugacy classes of $G$, with $C_0=\{1\}$.

It is proved in \cite{DE} that
$$\frac{A}{[A,A]}=HH_0(A)=\(HH_0\(\C[V], \C[V]\#G\)\)^G,$$
therefore
$$\frac A{[A,A]}=\(\bigoplus_{g\in G}HH_0\(\C[V], \C[V]g\)\)^G=\(\bigoplus_{g\in G}\C[V^g]g\)^G=\bigoplus_i\C[V^{g_i}]^{G_{g_i}}g_i$$
where $g_i$ is an arbitrary element of $C_i$ and $G_g=\{h\in G|gh=hg\}.$  The first summand is precisely $\C[V]^{G}$. Let $P_i$ be the projection $A\to \C[V^{g_i}]^{G_{g_i}}g_i$.

The Poisson bracket gives us a family of derivations $d_{\psi}:\C[V]^G\to \C[V]^G, \phi\mapsto\{\psi, \phi\}$ for $\psi\in\C[V]^G$; and we want to extend it to a larger family. The following Lemma restricts the possibilities.

\begin{lemma}
Let $d:A\to A$ be any derivation. If $x\in\C[V]^g\neq \C[V]$ then
$\(d(x)\)_g=0.$
\end{lemma}

\begin{proof}
Let $y\notin\C[V]^g.$ $d(xy)=d(yx)$ implies
$d(x)y+xd(y)=d(y)x+yd(x)$. The $g$-part of this equality is
$$\(d(x)\)_ggy+x\(d(y)\)_gg=\(d(y)\)_ggx+y\(d(x)\)_gg$$
or
$$\(d(x)\)_g \(^gy\)+x\(d(y)\)_g=\(d(y)\)_g \(^gx\)+y\(d(x)\)_g.$$
Since $^gx=x, ^gy\neq y$ we conclude $\(d(x)\)_g \(^gy-y\)=0$, so
$\(d(x)\)_g=0$
\end{proof}

Therefore if the action of $G$ on $V$ is faithful and $g\neq1$, the
$g$-part of the derivative an element of $\C[V]^g$ is zero. This
implies that for every $\psi\in \C[V]^G$, $d(\psi)\in \C[V]\subset
A$.

The condition $\langle
\overline{\psi g},\overline{\phi h} \rangle=-\langle \overline{\phi
h},\overline{\psi g} \rangle$ implies $\overline{d_{\psi g}(\phi
h)}=-\overline{d_{\phi h}(\psi g)}$. Consider the case $\phi, \psi\in\C[V]^G,\ h=1$ and $g\in C_i, i\neq0$. Since $P_i\(d_{\psi g}(\phi)\)=0$, we must have $0=P_i\(d_{\phi}(\psi
g)\)=P_i\(d_{\phi}(\psi) g+\psi d_{\phi}(g)\)$. The only terms that must be taken into account are
$\displaystyle d_{\phi}(\psi) g+\psi \sum
\(d_{\phi}(g)\)_{hgh^{-1}}hgh^{-1}$. Modulo $[A,A]$ this is equal to
\begin{multline*}
\(d_{\phi}(\psi) +\sum_h\ ^h\(\psi
\(d_{\phi}(g)\)_{hgh^{-1}}\)\)g=\(d_{\phi}(\psi) +\sum_h\psi\ ^h\(
\(d_{\phi}(g)\)_{hgh^{-1}}\)\)g\\=\(d_{\phi}(\psi) +\psi\sigma_{\phi,g} \)g
\end{multline*}
where $\displaystyle\sigma_{\phi,g}=\sum_h\ ^h\(
\(d_{\phi}(g)\)_{hgh^{-1}}\)$ does not depend on $\psi.$

We want $0=P_i\(\(d_{\phi}(\psi) +\psi\sigma_{\phi, g}
\)g\)=P_i\(\(\{\phi,\psi\} +\psi\sigma_{\phi,g} \)g\)$ since we want a
Poisson structure extending the usual Poisson bracket on $\C[V]^G.$ Therefore a neccesary condition for the existance of the Poisson structure is the existance of $\sigma_{\phi,g}\in\C[V]$ so that
\begin{equation}\label{condition}
P_i\(\(\{\phi,\psi\} +\psi\sigma_{\phi,g} \)g\)=0
\end{equation}
for every $\psi\in\C[V]$. We well see that this is not always possible.

\section{The counterexample}

Let $V=\C^4$ with linear coordinates $\{x_1, x_2, x_3, x_4\}$ and
the symplectic form $\omega=dx_1\wedge dx_2+dx_3\wedge dx_4$, so
$\C[V]=\C[x_1,x_2, x_3, x_4]$, and
$$\{\phi, \psi\}=\frac{\partial \phi}{\partial x_1}\frac{\partial \psi}{\partial x_2}-\frac{\partial \phi}{\partial x_2}\frac{\partial \psi}{\partial x_1}+\frac{\partial \phi}{\partial x_3}\frac{\partial \psi}{\partial x_4}-\frac{\partial \phi}{\partial x_4}\frac{\partial \psi}{\partial x_3}.$$

Let $G=\mathbb Z_2\ltimes\(\mathbb Z_2\oplus\mathbb Z_2\)$. Let $e,
b, c$ be the generators of the three copies of $\mathbb Z_2$ (in
that order). $G$ acts on $V$ as follows: $b$ and $c$ act as
$diag(-1, -1, 1, 1)$ and $diag(1, 1, -1, -1)$, respectively, on
$\{x_1, x_2, x_3, x_4\}$  and $e$ interchanges $x_1\leftrightarrow
x_3, x_2\leftrightarrow x_4$. Using Magma
(http://magma.maths.usyd.edu.au) we find that the ring of invariant
polynomials is generated, as an algebra, by $f_1=x_1^2 + x_3^2,\
f_2=x_2^2 + x_4^2,\ f_3=x_1^4 + x_3^4,\ f_4=x_2^4 + x_4^4,\
h_1=x_1x_2 + x_3x_4,\ h_2=x_1^2x_2^2 + x_3^2x_4^2,\ h_3=x_1^2x_3x_4
+ x_1x_2x_3^2,\ h_4=x_1x_2x_4^2 + x_2^2x_3x_4;$ with relations

$$-f_1f_2h_1 + f_1h_4 + f_2h_3 - h_1^3 + 2h_1h_2,$$
$$1/2f_1^2f_2 + 1/2f_1h_1^2 - 1/2f_1h_2 - 1/2f_2f_3 - h_1h_3,$$
$$1/2f_1f_2^2 - 1/2f_1f_4 + 1/2f_2h_1^2 - 1/2f_2h_2 - h_1h_4,$$
$$-1/2f_1^2f_4 + f_1f_2h_2 - 1/2f_2^2f_3 + f_3f_4 - h_2^2,$$
$$-1/2f_1^2h_4 + 1/2f_1f_2h_3 + 1/2f_1h_1h_2 - 1/2f_2f_3h_1 + f_3h_4 -h_2h_3,$$
$$1/2f_1f_2*h_4 - 1/2f_1f_4h_1 - 1/2f_2^2h_3 + 1/2f_2h_1h_2 + f_4h_3 -h_2h_4,$$
$$1/2f_1^3f_2 + 1/2f_1^2h_1^2 - f_1^2h_2 - 1/2f_1f_2f_3 - 1/2f_3h_1^2 +f_3h_2 -h_3^2,$$
$$1/2f_1^2f_2^2 - 3/4f_1^2f_4 + 1/2f_1f_2h_1^2 - 3/4f_2^2f_3 + f_3f_4 -1/2h_1^2h_2 - h_3h_4,$$
$$1/2f_1f_2^3 - 1/2f_1f_2f_4 + 1/2f_2^2h_1^2 - f_2^2h_2 - 1/2f_4h_1^2 +f_4h_2 -h_4^2.$$

\begin{prop}
The Poisson bracket on $\C[V]^G$ cannot be extended to a Poisson structure on $\C[V]\#G$ for $V$ and $G$ as defined above.
\end{prop}

\begin{proof}
Take $\phi=x_1^2+x_3^2,
\psi=x_1x_2+x_3x_4\in \C[V]^G$ and $g=b.$ In this case $C_i=\{b,c\}$
and $\C[V^b]^{G_b}b=\C[x_3^2, x_4^2, x_3x_4]b.$
$\{\phi,\psi\}=2x_1^2+2x_3^2$ so
$P_i\(\(\{\phi,\psi\}\)b\)=2x_3^2b$. On the other hand,
$P_i\(\(\psi\sigma_{\phi,g} \)b\)=P_i\(\(\(x_1x_2+x_3x_4\)\sigma_{\phi,g}
\)b\)=$\newline$P_i\(\(\(x_3x_4\)\sigma_{\phi,g} \)b\)$ and none of the terms here can
be equal to $-2x_3^2$ since they all contain $x_4.$ This contradicts (\ref{condition}).
\end{proof}

\section{Aknowledgments}
The author is grateful to William Crawley-Boevey for a careful reding of this paper, and to Victor Ginzburg for posing the question.

\end{document}